\documentclass[a4paper,11pt]{article}

\usepackage{mathtools}
\usepackage{enumerate}
\usepackage{verbatim}
\usepackage{mathrsfs}
\usepackage{float}
\usepackage[dvips]{graphicx}
\usepackage[latin1]{inputenc}
\usepackage{amsfonts}
\usepackage{amsmath}
\usepackage{amssymb}
\usepackage{amsthm}
\usepackage{epsfig} 
\usepackage[dvipsnames, usenames]{color}
\usepackage{geometry}
\usepackage{graphicx}
\usepackage[all]{xy}
\usepackage{cite}
\usepackage[usenames,dvipsnames,svgnames,table]{xcolor}
\usepackage{bbm}
\usepackage{leftidx}
\usepackage[normalem]{ulem}
\usepackage{algorithmic,algorithm}
\usepackage{hyperref}
\usepackage{titlesec} 

\DeclareMathOperator*{\argmax}{argmax}
\DeclareMathOperator{\sgn}{sgn}
\newcommand{\defeq}{\vcentcolon=}
\newcommand{\backdefeq}{=\vcentcolon}
\allowdisplaybreaks

\newtheorem{Theorem}{Theorem}[subsection]

\newtheorem{Assumptions}{Assumptions}[subsection]

\theoremstyle{plain}
\allowdisplaybreaks

\date{}
\begin{document}
\title{A policy iteration algorithm for nonzero-sum stochastic impulse games}

\author{Ren\'e Aïd\footnote{Université Paris Dauphine, PSL Research University, LEDA.}
	\qquad Francisco Bernal\footnote{CMAP, Ecole Polytechnique.} 
	\qquad Mohamed Mnif\footnote{University of Tunis El Manar, ENIT-LAMSIN.}\\
	Diego Zabaljauregui\footnote{Corresponding Author. Department of Statistics, London School of Economics and Political Science.\newline E-mail: {\tt d.zabaljauregui@lse.ac.uk}}
	\qquad Jorge P. Zubelli\footnote{IMPA.}}

\maketitle

\begin{abstract}
\noindent This work presents a novel policy iteration algorithm to tackle nonzero-sum stochastic impulse games arising naturally in many applications. Despite the obvious impact of solving such problems, there are no suitable numerical methods available, to the best of our knowledge. Our method relies on the recently introduced characterisation of the value functions and Nash equilibrium via a system of quasi-variational inequalities. While our algorithm is heuristic and we do not provide a convergence analysis,	numerical tests show that it performs convincingly in a wide range of situations,	including the only analytically solvable example available in the literature at the time of writing.
\end{abstract}	

{\em \qquad Accepted for publication in ESAIM: Proceedings and Surveys (November 2018)}\newline

\noindent {\em Keywords:} stochastic impulse game, nonzero-sum game, Nash equilibrium, policy iteration, Howard's algorithm, quasi-variational inequality.

\section*{Introduction}
\addcontentsline{toc}{section}{Introduction}

Stochastic impulse games (SIGs) are at the intersection between differential game theory and stochastic impulse control. In the case of one sole player, they reduce to stochastic impulse optimization problems where an agent seeks to control an {\em underlying} (or {\em state variable})---otherwise governed by a stochastic differential equation---in order to maximise the expected value of some target functional. One way in which she can influence the underlying is this: whenever it leaves a {\em continuation region} (i.e. at her {\em intervention times}), she suddenly shifts it back somewhere into the region (by providing an {\em impulse}). Together, the continuation region and impulses define what we will call a {\em strategy}. Thus, the control-theoretical problem typically boils down to determining an optimal strategy as a function of the current value of the underlying. The purpose of the analysis is to apply that optimal strategy to any specific realization of the underlying. When such a strategy is followed, the expected value of the target functional---the largest possible---is called the {\em value function}. This problem is a classical one and it is well understood by now \cite{OS}. In particular, a rigorous framework has been put in place---combining Howard's algorithm \cite{Be1,Be2,H,BMZ} with viscosity-capturing finite difference schemes---which allows for robust numerical approximations \cite{CMS,COS,OS}. 

A natural extension are two-player SIGs. In such games, two players seek to control the underlying with different, and often opposed, aims. Let us give two concrete examples, due to \cite{ABCCV_V2}.
\begin{itemize}
\item Two central banks competing to influence the exchange rate between their respective currencies, both seeking to devalue their own one. In this application, the exchange rate is modelled as a stochastic process, and either central bank intervenes when it deems its currency too strong. Each bank's strategy is made up of the ensuing devaluation and the threshold exchange rate triggering it. Both quantities are to be optimally determined in advance, by solving the game.  

\item A wholesale producer of energy and a big client of hers who resells it. Due to the high consumption of the client, both players are capable of affecting the wholesale price of energy, with naturally opposed targets (a high price versus a low one, respectively). Each player incurs costs when modifying the price while possibly benefiting from state compensation upon adverse price movements led by the opponent. In order  to avoid superfluous costs, both players need a strategy, based on the current wholesale price, determining the optimal triggering threshold and the extent of the price shift. 
\end{itemize}

Contrary to one-player SIGs, the value functions in two-player SIGs (one for each player) cannot be naively defined by maximisation. 
Instead, optimality can naturally be characterized via the notion of {\em Nash equilibrium}: intuitively, the pair of value functions displays the best expected outcome in the sense that if a player changes her strategy while her opponent does not, then the former can only be worse off. A pair of strategies at which these values are attained is then called a Nash equilibrium. Because of this, two-player SIGs are distinctly different (and more complex) from one-player SIGs, and the theoretical/numerical frameworks for the latter are not suitable for the former. This is further accentuated in the so-called \textit{nonzero-sum} SIGs (NZSSIGs), roughly meaning that the losses (resp. gains) undergone by one player may not exactly translate into gains (resp. losses) for the opponent. (For instance, the previously described examples are best modelled under this very flexible framework.) NZSSIGs are in contraposition with the simpler and more particular \textit{zero-sum} SIGs, in which both value functions add up to zero---effectively reducing the problem to finding one. (This allows for more natural extensions of the theoretical and numerical tools of the one-player case \cite{C,B}.)

Due to the greater complication involved in the analysis, NZSSIGs are much underdeveloped. Very recently, a breakthrough has been achieved in \cite{ABCCV_V2}, where the value functions and a Nash equilibrium of a very general class of NZSSIGs (assuming they exist and possess some regularity) have been characterized via a system of {\em quasi-variational inequalities} (QVIs). Besides the theoretical interest of this link by itself, it opens up the possibility of finding approximations to the value functions by numerically solving that system.

In this paper, we make the first attempt (to the best of our knowledge) at numerically solving NZSSIGs. In a nutshell, the iterative algorithm we put forward treats the NZSSIG at each iteration as a combination of a fixed point problem and a slowly relaxing one-player SIG. This allows us to take advantage of the machinery for the latter. Because, at the moment, we lack a proof of convergence, the proposed algorithm is admittedly heuristic. Instead, we report on a range of numerical experiments, which show convergence of the error with respect to the discretisation. In fact, errors are quite satisfactory and---in the examples we have tackled---the relative errors easily drop well below $0.1\%$. The algorithm can thus be used to assist further development of the field, as well as to gain insight into applications modelled by NZSSIGs. For even crafted NZSSIGs with exact solution are hard to construct: in fact, the authors of \cite{ABCCV_V2} also provide (as far as we know) the first example in the literature---thanks to which we were able to validate our algorithm in the first place. (For the sake of completeness, this solution has been included in Section \ref{Benchmark}.) Moreover, this example also evidences that even when an analytical solution is available, its computation in practice may involve parameters which require solving complex nonlinear systems of equations---thus stressing the convenience of numerical approximations.

Let us outline the organization of the rest of the paper. We start in Section 2 by properly setting up NZSSIGs and recalling the main result of \cite{ABCCV_V2}, namely, the Verification Theorem with the corresponding system of QVIs. For the sake of illustration, the concrete application of competition in energy markets is presented according to this framework. In Section 3, we briefly review the state-of-the-art numerics for one-player SIGs, with an emphasis on Howard's algorithm, which will later become a pivotal ingredient of our own numerical method for NZSSIGs. The latter is motivated, listed, and discussed in Section 4. Section 5 presents numerical evidence supporting it, in the light of which we draw some conclusions in Section 6.  

\section{ The theoretical model}

\subsection{ Two-player nonzero-sum stochastic impulse games}
\label{NZSSIG}
In this section we introduce a general class of two-player NZSSIG, to which our numerical algorithm is applicable. For the sake of briefness and pertinence, we will skip some of the most technical matters of the rigorous construction of the model. We refer the interested reader to \cite{ABCCV_V2} for more details and a more general modelling framework. 

Let $(\Omega, \mathcal F, (\mathcal F_t)_{t\geq 0}, \mathbb P)$ be a filtered probability space under the usual conditions, that supports a standard one-dimensional Brownian motion $W=(W_t)$. For each $x\in\mathbb R$ we consider a process $X$ (the \textit{state variable}) starting at $x$, 
\begin{equation}
\label{X}
X_t=x+ \int_0^t\mu(X_s)ds + \int_0^t\sigma(X_s)dW_s + \displaystyle\sum_{k:\ \tau^1_k\leq t}\delta^1_k + \displaystyle\sum_{k:\ \tau^2_k\leq t}\delta^2_k,
\end{equation}
where:
\begin{enumerate}[(i)]
\item $\mu,\sigma:\mathbb R\to \mathbb R$ are Lipschitz continuous.
\item $u^i\defeq\{(\tau^i_k,\delta^i_k)\}_{k=1}^\infty$ ($i=1,2$) is the \textit{control} of player $i$, which consists of a sequence of stopping times $\tau^i_k$ (the $k$-th \textit{intervention time} of player $i$) and $\mathcal F_{\tau^i_k}$-measurable random variables $\delta^i_k$ (the $k$-th \textit{intervention impulse} of player $i$).
\end{enumerate}
In words, when none of the players intervenes the state variable behaves as an It\^o diffusion. The players can at any time decide to shift this process by applying a certain impulse. The specific type of interventions are determined by \textit{(threshold-type) strategies}, which means that each player acts only when the state variable exits a given region of $\mathbb R$. More specifically, the control of player $i$ is defined in terms of a strategy $\varphi_i\defeq (\mathcal C_i,\xi_i),$ where $\mathcal C_i\subseteq\mathbb R$ is an open set (the \textit{continuation region}) and $\xi_i:{\mathcal C_i}^c\to \mathbb R$ is a continuous function.\footnote{For any $A\subseteq\mathbb R$, $A^c$ denotes the complement of $A$ in $\mathbb R$.} Player $i$ intervenes if $X$ exits $C_i$---i.e., if for some $t\geq 0$ and $\omega\in\Omega$ it holds $X_t(w)\notin\mathcal C_i$---by applying an impulse $\xi_i(X_t(w))$. We decree that the game never stops and if both players want to intervene at the same time, then player $1$ has the priority. The latter assumption is a matter of convention and not very restrictive.\footnote{Under the assumptions of the Verification Theorem \ref{verification}, the value of the game and the Nash equilibrium found would not change if player 2 had the priority instead. This can be checked by taking a pair $(\tilde V_1,\tilde V_2)$ as in the theorem and noticing that $(\tilde V_2,\tilde V_1)$ solves the game in which the structures (costs, gains, rates and payoffs) of the players have been swapped without swapping the priority.}

Henceforth, we write $\mathbb E_x$ to denote the conditional expectation given $X_{0^-}=x$, and define the \textit{objective function} for player $i$ given the strategies $(\varphi_1,\varphi_2)$ and the starting point $x$ of $X$ as
\begin{equation}
\label{Vi_given_strats}
\begin{split}
J^i(x;\varphi_1,\varphi_2)& \defeq \mathbb E_x \left[ \int_0^\infty e^{-\rho_i s} f_i(X_s)ds + \sum_{k=1}^\infty e^{-\rho_i \tau^i_k} \phi_{i} \big( X_{(\tau^i_k)^-}, \delta^i_k \big) + \sum_{k=1}^\infty e^{-\rho_i \tau^j_k} \psi_{i} \big( X_{(\tau^j_k)^-}, \delta^j_k \big)\right],
\end{split}
\end{equation}
with $i,j\in\{1,2\},\ j \neq i$; where:
\begin{enumerate}[(i)]
\item $\rho_i>0$ is the subjective \textit{discount rate} of player $i$.
\item $f_i:\mathbb R\to\mathbb R$ is a continuous function giving the \textit{running payoff} of player $i$.
\item $\phi_i:\mathbb R^2\to\mathbb R$ (resp. $\psi_{i}:\mathbb R^2\to\mathbb R$) is a continuous function giving the \textit{cost of intervention} (resp. \textit{gain due to the opponent's intervention}) for player $i$.
\end{enumerate}
In addition, we will only consider pairs of strategies $(\varphi_1,\varphi_2)$ such that the previous expectations are well defined for all $x\in\mathbb R$, and we refer to these pairs as \textit{admissible} (see \cite{ABCCV_V2} for more details).\footnote{In \cite{ABCCV_V2} admissibility is defined pointwise for greater generality. We refrain from doing this for simplicity.} Note for example that, for payoffs with polynomial growth, the `no intervention strategies' $\phi_1=\phi_2=(\mathbb R,\emptyset)$ are always admissible. In this context, we say that the players behave optimally if their strategies form a Nash equilibrium. We recall that $(\varphi_1^*,\varphi_2^*)$ is a \textit{Nash equilibrium} if it is an (admissible) pair of strategies such that, for all $(\varphi_1,\varphi_2)$,
$$
J^1(x;\varphi_1^*,\varphi_2^*)\geq J^1(x;\varphi_1,\varphi_2^*)\quad\mbox{and}\quad J^2(x;\varphi_1^*,\varphi_2^*)\geq J^2(x;\varphi_1^*,\varphi_2),
$$
i.e., if player $i$ changes strategy while player $j$ does not, then on average player $i$ will be worse off, and vice versa. If a Nash equilibrium $(\varphi_1^*,\varphi_2^*)$ exists, we define the \textit{value function} for player $i$ when using the strategies $(\varphi_1^*,\varphi_2^*)$ by 
\begin{equation}
\label{Vi}
V_i(x)\defeq J^i(x;\varphi_1^*,\varphi_2^*),\qquad i=1,2.
\end{equation}
Our aim is to compute $(V_1,V_2)$ for some Nash equilibrium, and more importantly, to retrieve from these values the equilibrium itself.

\subsection{ The system of quasi-variational inequalities}

In order to establish a system of QVIs for $(V_1,V_2)$ we need to define one last ingredient, known as the \textit{intervention operators}, which will display the effect of each player's intervention on the value functions. For any two arbitrary functions $\tilde V_1, \tilde V_2:\mathbb R\to\mathbb R$, $i,j \in \{1,2\},\ i \neq j$, and $x\in\mathbb R$, the \textit{loss operator} of player $i$ is given by 
\begin{equation}
\label{Mi}
\mathcal M_i\tilde V_i(x)\defeq\displaystyle\sup_{\delta\in\mathbb R}\{\tilde V_i(x+\delta) + \phi_i(x,\delta)\}.
\end{equation}
If $\tilde V_i=V_i$ this operator gives the recomputed present value of $i$ due to the cost of her own intervention. 
If for each $x\in\mathbb R$ there exists a unique $\delta^j(x)=\delta^j_{\tilde V_j}(x)$ that realizes the supremum in (\ref{Mi})
, we also define the \textit{gain operator} of player $i$ as
\begin{equation}
\label{Hi}
\mathcal H_i\tilde V_i(x)\defeq \tilde V_i(x+\delta^j(x)) + \psi_i(x,\delta^j(x)).
\end{equation}
If $\tilde V_i=V_i$ this operator gives the recomputed present value of player $i$ due to her opponent's intervention. Whenever we make use of a gain operator, we are implicitly stating that the above assumptions and notations are in place. We also emphasize that $\mathcal H_i\tilde V_i(x)$ depends on the whole function $\tilde V_j$ through $\delta^j_{\tilde V_j}$ and we will write $\mathcal H_i(\tilde V_j)\tilde V_i(x)$ instead, when we want to make this dependence explicit.

We can now state the Verification Theorem, due to \cite{ABCCV_V2}, that will allow us to tackle the problem of finding $(V_1,V_2)$ and a Nash equilibrium by numerically solving a deterministic system of QVIs. In the next theorem we use the notation $\mathcal A$ for the infinitesimal generator of $X$ when no interventions take place. That is, $\mathcal A$ is the operator such that if $g\in C^2(S)$ for some $S\subseteq\mathbb R$, then
\begin{equation}
\label{A}
\mathcal Ag(x)= \mu(x)g'(x) + \frac{1}{2}\sigma^2(x) g''(x),\quad\mbox{ for }x\in S.
\end{equation}
We concur that whenever we apply $\mathcal A$ to some function $g$, we are implicitly stating this function $g$ is $C^2$ at every $x$ at which we compute $\mathcal Ag(x)$.

\begin{Theorem}[\textbf{System of QVIs}]
\label{verification}
Let $\tilde{V}_1,\tilde{V}_2:\mathbb R\to\mathbb R$ such that for any $i,j \in \{1,2\}$, $i \neq j$:
\begin{equation}
\label{QVIs}
\begin{cases}
	\begin{aligned}
		& \mathcal M_j\tilde V_j -\tilde V_j \leq 0, && \text{in} \,\,\, \mathbb R,   \\
		& \mathcal H_i\tilde V_i-\tilde V_i=0, && \text{in} \,\,\, \{\mathcal M_j\tilde V_j - \tilde V_j = 0\},  \\
		& \max\big\{\mathcal A \tilde V_i -\rho_i \tilde V_i + f_i, \mathcal M_i\tilde V_i-\tilde V_i \}=0, && \text{in} \,\,\, \{\mathcal M_j\tilde V_j -\tilde V_j < 0\}\backdefeq\mathcal C_j,
	\end{aligned}
\end{cases}
\end{equation}
and $\tilde{V}_i\in C^2(\mathcal C_j\backslash\partial\mathcal C_i)\cap C^1(\mathcal C_j)\cap C(\mathbb R)$ has polynomial growth and bounded second derivative on some reduced neighbourhood of $\partial\mathcal C_i$. Suppose further that $(\varphi_1^*,\varphi_2^*)$, with $\varphi_i^*\defeq(\mathcal C_i,\delta^i_{\tilde V_i})$, is an admissible pair of strategies. Then 
$$(\tilde{V}_1,\tilde{V}_2)=(V_1,V_2)\mbox{ and }(\varphi_1^*,\varphi_2^*)\mbox{ is a Nash equilibrium.}$$
\end{Theorem}
\subsection{ Example of application: competition in energy markets}\label{Benchmark}
We finish this Section by providing a specific example of application: competition in energy markets \cite{ABCCV_V2}.

Let process $X$ model the forward price of energy, evolving as a Brownian motion when there are no interventions. Consider the following two players: player 1 is an energy producer with unitary production cost $s_1$, so that in a simplified model her payoff is $X-s_1$. Player 2 runs a large company that buys from player 1 and sells at a unitary price $s_2>s_1$, with payoff $s_2-X$. Because of the high consumption of player 2, she can also affect the price of energy in the same way as player 1. Whenever these players intervene they incur a cost (advertising among other factors) which is modelled linearly for tractability---one constant component and another proportional to the change induced in the energy price. At the same time, because of their impact in the economy, the government subsidizes upon adverse movements in the energy price. For each player this represents a gain when the opponent intervenes, which is modelled in the same way as the cost of intervention, but with different parameters. We further assume that both players discount their winnings/losses at the same rate $\rho>0$ and they have the same cost/gain parameters. More specifically,
$$
X_t=x+ \sigma W_t + \displaystyle\sum_{k:\ \tau^1_k\leq t}\delta^1_k + \displaystyle\sum_{k:\ \tau^2\leq t}\delta^2_k
$$
and for $i,j\in\{1,2\},\ i \neq j$,
\begin{equation}
\label{J_retail}
J^i(x;\varphi_1,\varphi_2) \defeq \mathbb E_x \left[ \int_0^\infty e^{-\rho s} (-1)^{i-1}(X_s - s_i)ds - \sum_{k=1}^\infty e^{-\rho \tau^i_k} (c+\lambda |\delta^i_k|)+ \sum_{k=1}^\infty e^{-\rho \tau^j_k} (\tilde{c}+\tilde{\lambda} |\delta^j_k|)\right],
\end{equation}
where $0\leq \tilde{c}\leq c,\ 0\leq\tilde{\lambda}\leq\lambda,\ (c,\lambda)\neq (\tilde{c},\tilde{\lambda}),\mbox{ and } 1-\rho\lambda>0$. These parametric restrictions ensure among other things the existence of a Nash equilibrium (see \cite{ABCCV_V2} for more details). 

The exact solution to this game can be found by application of Theorem \ref{verification}. To this purpose, the system of QVIs (\ref{QVIs}) is heuristically solved, yielding candidates for value functions and a Nash equilibrium. This is done, first, by making some educated guesses regarding the shape of the optimal continuation regions and value functions. Second, by solving the ordinary differential equations (ODEs) in (\ref{QVIs}) where appropriate. Finally, by imposing the regularity requirements of Theorem \ref{verification} through pasting conditions. Upon verification of the remaining hypotheses, the following turns out to be the solution to the game:
\begin{align*}
&V_2(x)= \left\{
\begin{array}{ll}
\varphi^{A_{1},A_{2}}(x_1^*) + {\tilde c} + {\tilde{\lambda}}(x_1^*-x) &\textrm{ if } x\in(-\infty,{\bar x}_1],\\
\varphi^{A_{1},A_{2}}(x) &\textrm{ if } x\in({\bar x}_1,{\bar x}_2),\\
\varphi^{A_{1},A_{2}}(x_2^*) - c - \lambda(x-x_2^*) &\textrm{ if } x\in[{\bar x}_2,+\infty),
\end{array}
\right.
&V_1(x)= V_2(2{\tilde s}-x),	
\end{align*}
where:
\begin{align*}
&\varphi^{A_{1},A_{2}}(x)= A_{1}e^{\theta x} + A_{2}e^{-\theta x} + \frac{1}{\rho}(s_2-x),\\
&{\tilde s}:= \frac{s_1+s_2}{2}, \qquad
\theta:= \sqrt{2\rho/\sigma^2},\qquad
\eta:= (1-\lambda\rho)/\rho,\\
&{\bar x}_i:= {\tilde s} + \frac{(-1)^i}{\theta}\log{\Bigg(\sqrt{\frac{\eta+\xi}{\eta-\xi}}\Big(\sqrt{\Gamma+1}+\sqrt{\Gamma}\Big)\Bigg)},\quad x^*_i:= {\tilde s} + \frac{(-1)^i}{\theta}\log{\Bigg(\sqrt{\frac{\eta-\xi}{\eta+\xi}}\Big(\sqrt{\Gamma+1}+\sqrt{\Gamma}\Big)\Bigg)},\\
&A_{i}:=\exp{\Big((-1)^{i}\theta {\tilde s}\Big)}
\frac{\sqrt{\eta^2-\xi^2}}{2\theta}
\Big((-1)^{i+1}\sqrt{\Gamma+1}-\sqrt{\Gamma}\Big),\\
&\Gamma:= \frac{\theta(c-{\tilde c})}{4\xi}+
\frac{\theta c(\lambda-{\tilde{\lambda}})}{4\eta\xi}+
\frac{\lambda-{\tilde\lambda}}{2\eta}
\end{align*}
and $\xi\in(0,\eta)$ is the unique zero of $F(y):= 2y  - \eta\log{\Big(\frac{\eta+y}{\eta-y}\Big)} + \theta c$. 
\bigskip

As it can be readily noticed, the analytical solution involves the computation of several parameters and the resolution of at least one nonlinear equation. As a matter of fact, the number of parameters in this solution was reduced making use of the symmetry in the problem. In a more general two-player NZSSIG under the modelling framework of Section \ref{NZSSIG}, if we assume the optimal continuation regions in Theorem \ref{verification} are semi-bounded intervals, then an analytical solution can easily involve eight parameters: two finite limits of both continuation regions, two maximising points of the net value functions (subtracting intervention costs) and four undetermined constants coming from the second order ODEs. Moreover, these parameters are the solution to a nonlinear system of equations that arises from smooth pasting and optimality conditions \cite[Def. 4.1.]{ABCCV_V2}. In the more frequent than not situation in which this system cannot be simplified, computing the analytical solution may become prohibitive. This further motivates the need for a numerical algorithm to solve the system of QVIs (\ref{QVIs}). \newline

\noindent
{\bf Remark.} To the best of our knowledge this is, at the time of writing, the only one example in the literature of an analytically solvable NZSSIG. Henceforth, we will refer to it as the {\em benchmark game}.

\section{Numerics for one-player SIGs: state of the art}
\setcounter{subsection}{1}
\noindent
{\bf Notation.} Henceforward, finite grids will be denoted by $\mathbb K$ and $\mathbb S$. For a fixed grid, we will use the same font type when discretising an operator, e.g., $\mathcal M,\ \mathcal H\mbox{ and }\mathcal A$ will be replaced by $\mathbb M,\ \mathbb H\mbox{ and } \mathbb A$. We shall specify later the way the discretisations have been done in our experiments. We will not change notation for functions over grids as their domain will always be clearly stated. The only exception will be those functions that have been redefined to account for boundary conditions (BCs). Lastly, we recall that for any set $A$, $\mathbb R^A$ denotes the set of functions from $A$ to $\mathbb R$.
\bigskip

\textit{Policy iteration} or \textit{Howard's} algorithm for (discrete) variational inequalities (VIs) was originally developed in \cite{Be1, Be2, H}. The method was then extended to QVIs \cite{CMS, COS, OS}, i.e., variational inequalities in which the obstacle depends on the solution itself. Let $\mathbb K\subset\mathbb R$ be a finite set (the grid). These problems have the form:  
\begin{equation}
\label{QVI}
\mbox{Find } V\in\mathbb R^{\mathbb K}\mbox{: }\quad\max\big\{\mathbb L V(x)+\mathbbm g(x),\ \max_{y\in\mathbb K}\mathbb B^y V(x)- V(x)\big\}=0\quad\mbox{for all }x\in\mathbb K,\footnote{Slightly more general formulations are also available.}
\end{equation}
where $\mathbbm g\in\mathbb R^{\mathbb K}$ and $\mathbb L,\mathbb B^y:\mathbb R^{\mathbb K}\to\mathbb R^{\mathbb K}$ are linear and affine operators resp., for each $y\in\mathbb K$. 

Problem \ref{QVI} encomprises in particular a discrete localized version of a one-player SIG (i.e., a stochastic impulse control problem) since the system (\ref{QVIs}) reduces in this case to one single QVI. Indeed, by an appropriate finite difference approximation (more details in Section \ref{section_our_algorithm}) one can take $\mathbbm g =f|_{\mathbb K}$ (modified to account for Dirichlet or Neumann BCs),\footnote{For any function $F:A\to B$ and $X\subseteq A$, $F|_X$ denotes the restriction of $F$ to $X$.} $\mathcal A - \rho_1 Id\approx\mathbb L$ and $\mathcal M_1\tilde V(x)\approx \max_{y\in\mathbb K}\mathbb B^y\tilde V(x)$, with $\mathbb B^y\tilde V(x)\defeq \tilde V(y)+\phi_1(x,y-x)$. By their antisymmetric structure, two-player zero-sum SIGs can also be accommodated with a similar formulation as a max-min double obstacle problem \cite[page 7]{ABCCV_V2} for one single value function $V\defeq V_1=-V_2$, and the policy iteration algorithm can be adapted \cite{BMZ}.

However, this framework is not general enough to tackle the NZSSIGs described in Section \ref{NZSSIG} and the full system of QVIs (\ref{QVIs}). There are, to the best of our knowledge, no available numerical methods to approach the latter problem. We present now the classical policy iteration algorithm (Algorithm \ref{Ho_QVI}) used to solve (\ref{QVI}). It will be embedded in the numerical scheme we will put forward to solve the much more general problem (\ref{QVIs}). While we note that our later use of this algorithm is not fully within the theoretical scope studied in \cite{CMS}---in particular, we will need to relinquish the affine nature of the operators $\mathbb B^y$, $y\in\mathbb K$---the two following main assumptions will be in place nonetheless ($\#{\mathbb K}$ is the cardinal of set ${\mathbb K}$):

\begin{Assumptions}
\label{assumptions}
\begin{enumerate}[(i)]
\item $-\mathbb L$ is a strictly diagonally dominant $M$-operator. That is, if $L\in\mathbb R^{\#\mathbb K\times\#\mathbb K}$ is the canonical matrix of $\mathbb L$,\footnote{If $\mathbb K$ is the grid: $x_0<x_1<\dots<x_M$, then $L$ is the matrix with columns $L_j=\mathbb L\mathbbm 1_{\{x_j\}}$, for $j=0,\dots,M$.} then 
$$L_{ij}\geq 0\mbox{ for all }i\neq j\mbox{ and }-L_{ii}>\sum_{j\neq i}L_{ij}\mbox{ for all }i.$$
\item $\mathbb B^y$ is a non-expansive function for $\|\cdot\|_\infty$, for all $y\in\mathbb K$. That is,
$$\|\mathbb B^y \tilde V_1 - \mathbb B^y \tilde V_2\|_\infty\leq \|\tilde V_1-\tilde V_2\|_\infty,\quad\mbox{ for all }\tilde V_1,\tilde V_2\in\mathbb R^{\mathbb K},\ y\in\mathbb K.$$
\end{enumerate}
\end{Assumptions}
\noindent We remark that, in particular, Assumption \textit{(i)} will guarantee that Algorithm \ref{Ho_QVI} is well defined, as the linear operators $\mathbb L^k$ will be non-singular \cite{BMZ}. 

Lastly, note that in Algorithm \ref{Ho_QVI}, computing the operators $\mathbb L^k$ amounts simply to redefining a matrix row by row, using either the matrix of $\mathbb L$ or $-Id$. We have chosen an `operators-type' notation for this paper, as it will simplify matters in the sequel, when compared with its matrix counterpart.

\begin{algorithm}[H]
\caption{Policy iteration for one QVI (one-player SIG)}
\begin{algorithmic}[1]
\STATE{Set $\varepsilon>0$ (numerical tolerance) and $k_{\max}\in\mathbb N$ (maximum iterations).}
\STATE{Pick initial guess: $V^0\in\mathbb R^\mathbb K$.}
\STATE{Let $k=0$ (iteration counter) and $R^0=+\infty$.}
\WHILE{$R^k>\varepsilon$ \AND $k\leq k_{\max}$} 
\STATE{$M^k \defeq \max_{y\in\mathbb K}\mathbb B^y V^k$.}
\STATE{$\alpha^k\defeq\mathbbm 1_{\{\mathbb L V^k+\mathbbm g <M^k-V^k\}}$ (action at each point).} 
\STATE{Define $\mathbb L^k:\mathbb R^{\mathbb K}\to\mathbb R^{\mathbb K}$ and $\mathbbm g^k\in\mathbb R^{\mathbb K}$ by\newline
\begin{equation*}
\mathbb L^k\tilde V(x)\defeq \begin{cases}
										\mathbb L\tilde V(x) & \mbox{ if }\ \alpha^k(x)=0\\
										-\tilde V(x) & \mbox{ if }\ \alpha^k(x)=1
										\end{cases}\quad
\quad
\mathbbm g^k(x)\defeq \begin{cases}
										\mathbbm g(x) & \mbox{ if }\ \alpha^k(x)=0\\
										M^k(x) & \mbox{ if }\ \alpha^k(x)=1.
									 \end{cases}
\end{equation*}
}
\STATE{Solve for $V^{k+1}$:  $\ \mathbb L^k V^{k+1} + \mathbbm g^k = 0$.}
\STATE{$R^{k+1}\defeq\|V^{k+1}-V^k\|$.}
\STATE{$k=k+1$.}
\ENDWHILE
\end{algorithmic}
\label{Ho_QVI}
\end{algorithm}

\section{ Proposed algorithm for two-player NZSSIGs} 
\label{section_our_algorithm}
\setcounter{subsection}{1}
Compared with the single-value-function problems in the previous Section, general two-player NZSSIGs are distinctly more challenging. The main challenges are:   
\begin{itemize}
	\item two value functions, governed by a system of QVIs, must be solved for,	
	\item the dependence between $(V_1,V_2)$ is highly nonlinear due to the presence of the gain operators $\mathcal H_i(V_j)V_i$,
	\item each gain operator is expansive as a function of $(\tilde V_1,\tilde V_2)$; and 
	\item solutions will typically be less regular. For example, if $(V_1,V_2)$ is the solution to the benchmark game in section \ref{Benchmark} then $V_i$ is singular at $\bar{x}_j$ for $i,j\in\{1,2\},\ i\neq j$ (i.e., each value function is non differentiable at the border of the opponent's continuation region), in spite of the game having linear payoffs, costs and gains. This is somehow to be expected, as Theorem \ref{verification} contemplates this lack of regularity within the smoothness assumptions (compare to the classical Verification Theorems for one-player problems \cite{OS} where greater regularity is assumed). 
\end{itemize}

Algorithm \ref{A:Our_algorithm} below is, as far as we know, the first ever numerical attempt at two-player NZSSIGs. At present, it is admittedly heuristic and supported only by the numerical evidence reported in Section \ref{S:Experiments}. 

The remainder of this Section is organized as follows. We start by explaining the idea and motivating the underlying heuristics. Then, we describe in detail the discretisation of the system of QVIs (\ref{QVIs}). Finally, we list the new algorithm.\newline

{\bf Heuristics.}
Having discretised the system of QVIs (\ref{QVIs}), we start from an initial guess $(V_1^0,V_2^0)$ to approximate its solution. We seek an iterative procedure to consecutively compute $(V_1^{k+1},V_2^{k+1})$, given  $(V_1^k,V_2^k)$ at the $k$-th iteration.

Let $i,j\in\{1,2\},\ i\neq j$. A natural idea is, first, to partition the grid into the `approximate continuation region' of player $j$---$\{\mathbb M_jV_j^k-V_j^k<0\}$---and its complement, the `approximate \textit{intervention region}'; and then to compute $V_i^{k+1}$ either by calculating a gain as $\mathbb H_i(V_j^k)V_i^k$ in the former case, or by solving one QVI with Howard's algorithm (Algorithm \ref{Ho_QVI}) in the latter. Note, however, that naively defining $\{\mathbb M_jV_j^k-V_j^k<0\}$ as the `approximate continuation region' of player $j$ poses difficulties. 

Indeed, the discretisation of the loss operator $\mathcal M_j$ as $\mathbb M_j$ implies that even the true value function $V_j$ will generally not verify $\mathbb M_jV_j(x)-V_j(x)=0$ (or $\mathbb M_jV_j(x)-V_j(x)\geq 0$) for $x$ in the true intervention region of player $j$. We therefore need to relax this constraint to account both for numerical error and the discrepancy between the discrete and space-continuous problems. Our experiments have shown that a successive relaxation procedure turns out to be the most effective. Consequently, we define the approximate continuation region of player $j$ as $\{\mathbb M_jV_j^k-V_j^k<-r^k\}$ instead, where $r^k$ is a small positive number. By letting $r^k$ relax to a preset small tolerance $\varepsilon>0$, the iterative approximations will hopefully converge to the correct discrete solution.

It remains to schedule the relaxation procedure and to define a measure of convergence for the algorithm. Regarding the former, having computed $(V_1^{k+1},V_2^{k+1})$, $r^k$ is linearly relaxed (line 11). (This relaxation procedure is chosen for simplicity.). Then the largest pointwise residual to the system of QVIs (incurred by either approximate value function) is calculated across the grid (line 12), taking into account the numerical tolerance $\varepsilon$. We denote this residual $R^{k+1}$ and we consider the algorithm has converged when it drops below a certain tolerance---unlike in Algorithm \ref{Ho_QVI}, where the residual is taken as the distance between consecutive approximations. This alternative approach has been chosen for being more informative. It reflects whether a solution to the discrete system of QVIs has been found---as opposed to the algorithm stagnating---on top of giving valuable information at each grid node. 

By construction, the junctions between the approximate continuation and intervention regions for each of the players will necessarily take place on top of a finite difference node. This may lead to numerical issues when the exact value functions are non differentiable there (as will often be the case). Namely, the pointwise residual to the QVIs at the junction nodes may not be made arbitrarily small by refining the grid---because the derivatives contained in the residual are not defined at the exact junction, in the first place. 
In those cases, the algorithm must be stopped at that point, since adding more grid nodes in a naive way cannot lead to any improvement. 

As a final remark, we mention that numerical experiments show Algorithm \ref{A:Our_algorithm} does not enjoy global convergence (i.e. it is not guaranteed to converge from arbitrary initial guesses). Providing a good enough pair $(V_1^0,V_2^0)$ is thus a practical issue; in Section \ref{S:Experiments}, a natural way of constructing educated guesses is explained.\newline 

\textbf{discretisation.}
Let $\mathbb S\subseteq\mathbb R$ be a finite set. $\mathbb S$ is the grid we will use to discretise the system (\ref{QVIs}). In all of the numerical experiments described in the sequel we have taken $\mathbb S$ as an equispaced grid of $M$ steps between certain $x_{\min}<0<x_{\max}$ with $|x_{\min}|,|x_{\max}|$ and $M$ big enough (more about this below).

Let $i,j\in\{1,2\},\ i\neq j$, $\tilde V_1,\tilde V_2\in\mathbb R^{\mathbb S}$ and $x\in\mathbb S$. We proceed to define the discretised versions of the loss and gain operators, over the grid $\mathbb S$, as  
$$
\mathbb M_i\tilde V_i(x)\defeq\max_{y\in\mathbb S}\{\tilde V_i(y)+\phi_i(x,y-x)\}\qquad\mathbb H_i(\tilde V_j)\tilde V_i(x)\defeq \tilde V_i(\mathbbm y_j(x))+\psi_i(x,\mathbbm y_j(x)-x),
$$
where 
$$
\mathbbm y_j(x)\defeq \min\left( \argmax_{y\in\mathbb S}\{\tilde V_j(y)+\phi_j(x,y-x)\}\right).
$$ 
Next, we choose a finite difference scheme for the ODEs
$$
0=\mathcal A \tilde V_i -\rho_i \tilde V_i + f_i =\frac{1}{2}\sigma^2 \tilde V_i''+\mu\tilde V_i'-\rho_i \tilde V_i + f_i
$$
which is consistent, monotone and stable, adding Dirichlet or Neumann-type BCs. In all of our experiments, we have chosen an upwind finite difference scheme, where
$$
\tilde V_i'(x)\approx \frac{\tilde V_i\big(x+\sgn(\mu(x)) h\big)-\tilde V_i(x)}{\sgn(\mu(x))h}\qquad \tilde V_i''(x)\approx \frac{\tilde V_i(x+ h) - 2\tilde V_i(x) + \tilde V_i(x-h)}{h^2},\footnote{$\sgn$ denotes the sign function, i.e., $\sgn(x)=1$ if $x\geq 0$ and $\sgn(x)=-1$ otherwise.}
$$
and $\tilde V_i(x_{\min}-h),\tilde V_i(x_{\max}+h)$ were solved for using Neumann conditions on 
$$
\tilde V_i'\left(x_{\min}-\frac{1+\sgn(\mu(x_{\min}))}{2}h\right),\quad\tilde V_i'\left(x_{\max}+\frac{1-\sgn(\mu(x_{\max}))}{2}h\right)
$$ 
respectively. How to get these conditions, and in particular how to choose $x_{\min}$ and $x_{\max}$, is a problem-specific question. In some situations, and particularly in the models numerically tested in this paper, one can heuristically assert that at a Nash equilibrium the continuation region of player $i$ should be a semi-interval of the form $\mathcal C_i=(\underline{x}_i,\overline{x}_i)$---with one end-point finite and the other infinite. Intuitively, one can further guess that there should exist a unique $y^*_i\in(\underline{x}_i,\overline{x}_i)$ that maximises the net value of player $i$ when she intervenes (see, e.g., \cite[Sect. 4.2]{ABCCV_V2} and \cite[Sect. 2.3.1]{B}). These conjectures arise mainly from the observation of the payoff functions $f_1, f_2$ which encode the goals of the players and roughly hint at some broad regions where each player would like the state variable to remain at. On $\mathbb R\backslash\mathcal C_i\mbox{ and }\mathbb R\backslash\mathcal C_j$
one of the two players will intervene, thus giving either 
$$V_i(x)=\mathcal M_iV_i(x) = V_i(y^*_i)+\phi_i(x,y^*_i-x)\quad\mbox{or}\quad V_i(x)=\mathcal H_iV_i(x)= V_i(y^*_j)+\psi_i(x,y^*_j-x).$$ 
If on the interiors $(\mathbb R\backslash\mathcal C_i)^o\mbox{ and }(\mathbb R\backslash\mathcal C_j)^o$ the derivatives $\frac{d\phi_i(x,y^*_i-x)}{dx}\mbox{ and } \frac{d\psi_i(x,y^*_j-x)}{dx}$ exist and do not depend on $y^*_i, y^*_j$ resp., differentiating the previous relations yields the Neumann conditions (provided $x_{\min},x_{\max}$ are `extreme' enough). The previous requirements on the derivatives are satisfied, for example, if the cost and gain structures have the form $\phi_i(x,\delta)=g_i(x)+ a_i|\delta|$ and $\psi_i(x,\delta)= h_i(x)+b_i|\delta|$, for some differentiable functions $g_i,h_i:\mathbb R\to\mathbb R$ and $a_i,b_i\in\mathbb R$. For the benchmark game
the Neumann BCs read $V_1'(x_{\min}-h)=\lambda,\ V_1'(x_{\max})=\tilde\lambda,\ V_2'(x_{\min}-h)=-\tilde\lambda\mbox{ and }V_2'(x_{\max})=-\lambda$. We will denote by $\mathbb A$ the discretised version of $\mathcal A$ and $\mathbbm f_i$ the restriction of $f_i$ to $\mathbb S$, redefined at $\min\mathbb S$ and $\max\mathbb S$ to account for the BCs. 

For Algorithm \ref{A:Our_algorithm} we recall that given any real number $a$, $a^+$ denotes its positive part (i.e., $a^+\defeq\max\{a,0\}$) and for any subset $S\subseteq \mathbb R$, $\mathbbm 1_{S}$ denotes the indicator function of $S$. 

\begin{algorithm}[H]
	\caption{Policy iteration for system of QVIs (two-player NZSSIG)}
	\begin{algorithmic}[1]
		\STATE{Set $\varepsilon>0$ (numerical tolerance), $0<\alpha<1$, $r^0>0$ (relaxation parameters) and $k_{\max}\in\mathbb N$ (maximum iterations).}
		\STATE{Pick initial guess: $(V_1^0,V_2^0)\in\mathbb R^{\mathbb S}\times \mathbb R^{\mathbb S}.$}
		\STATE{Let $k=0$ (iteration counter) and $R^0=+\infty$}
		\WHILE{$R^k>\varepsilon$ \AND $k\leq k_{\max}$} 
		\FOR{i=$1,2$ (player $i$)}
		\STATE{$j=3-i$ (player $j$).}
		\STATE{${\cal C}_j^k:=\{\mathbb M_jV_j^k-V_j^k<-r^k\}$.}
		\STATE{For $x\notin {\cal C}_j^k$, let $V_i^{k+1}(x)=\mathbb H_i(V_j^k)V_i^k(x)$.}
		\STATE{For $x\in {\cal C}_j^k$, solve for $V_i^{k+1}(x)$ by applying Algorithm \ref{Ho_QVI} to\newline
			$\max{\Big\{
				{\mathbb A}V_i^{k+1}(x)-\rho_i V_i^{k+1}(x)+\mathbbm f_i(x),\
				{\mathbb M}_iV_i^{k+1}(x)-V_i^{k+1}(x)	
				\Big\}=0}$.}
		\ENDFOR
		\STATE{$r^{k+1}\defeq\max{\{\alpha r^k,\varepsilon\}}$ (relaxation).}
		\STATE{Let $R^{k+1}$ be the largest pointwise residual to the system of QVIs, i.e.
			\begin{align}
			\label{eq:R^max}
			R^{k+1}:= \max_{\substack{i,j\in\{1,2\},j\neq i\\x\in \mathbb S}}&\Big\{
			\big({\mathbb M}_iV_i^{k+1}(x)-V_i^{k+1}(x)\big)^+,\\
			&\big|{\mathbb H}_i(V_j^{k+1})V_i^{k+1}(x)-V_i^{k+1}(x)\big|\mathbbm 1_{{\cal C}_j^{k,\varepsilon}}(x),\nonumber\\
			&\big|\max\big\{\mathbb A V_i^{k+1}(x)-\rho_i V_i^{k+1}(x)+\mathbbm f_i(x),{\mathbb M}_iV_i^{k+1}(x)-V_i^{k+1}(x)
			\big\}\big|\mathbbm 1_{\mathbb S\backslash{\cal C}_j^{k,\varepsilon}}(x)
			\Big\},\nonumber
			\end{align}		
			where ${\cal C}_j^{k,\varepsilon}:=\{ \mathbb M_jV_j^{k+1}-V_j^{k+1}<-\varepsilon\}$.}	
		\STATE{Let $k=k+1$.}
		\ENDWHILE	
	\end{algorithmic}
	\label{A:Our_algorithm}
\end{algorithm}	

Line 9 of Algorithm \ref{A:Our_algorithm} deserves some special attention. Although at this step we want to solve a problem restricted to the subgrid $\mathbb K={\cal C}_j^k$ (for fixed $k,i,j$), we still need the information in $\mathbb S\backslash {\cal C}_j^k$ in two ways, namely: 
\begin{itemize}
\item To compute the non-local operators $\mathbb M_i$. 
\item To restrict to $\mathbb K$ the equation $\mathbb A V_i^{k+1}(x)-\rho_i V_i^{k+1}(x)+\mathbbm f_i(x)=0$, properly accounting for the original BCs.
\end{itemize}
Suppose $\mathbb S$ is the grid $x_0<x_1<\dots< x_M$ and let us identify each point with its respective index. Let $A\in\mathbb R^{(M+1)\times(M+1)}$ be the canonical matrix of $\mathbb A$. For any subsets $I,J\subseteq\{0,\dots,M\}$, let us write $A_{I,J}$ for the submatrix of $A$ which has rows indexed in $I$ and columns indexed in $J$, and $\mathbb A_{I,J}$ the associated operator. 
Put $H^k_i\defeq V_i^{k+1}|_{\mathbb S\backslash \mathbb K}= \mathbb H_i(V_j^k)V_i^k|_{\mathbb S\backslash \mathbb K}$ and $h^k_i\defeq \max H^k_i.$  Then in order to apply Algorithm \ref{Ho_QVI} we take
$$
\mathbb L\tilde V=\mathbb A_{\mathbb K,\mathbb K}\tilde V-\rho_i \tilde V,\quad\mathbbm g=\mathbbm f_i|_\mathbb K + \mathbb A_{\mathbb K,\mathbb S\backslash\mathbb K}H^k_i\quad\mbox{and}\quad\mathbb B^y\tilde V(x) = \max\{\tilde V_i(y)+\phi_i(x,y-x),h^k_i\},
$$
for all $\tilde V\in\mathbb R^{\mathbb K},\ x,y\in\mathbb K$. We remark once again that the functions $\mathbb B^y$ fail to be affine operators. This, together with some assumptions which are not satisfied, make this application of Algorithm~\ref{Ho_QVI} fall outside of the scope of \cite{CMS} (and even more of \cite{COS, OS}). However, it is easy to check that the main assumptions, Assumptions \ref{assumptions}, are still verified under our discretisation and we have observed unconditional convergence of this subroutine (to very high orders of precision) in all the experiments we have ran. 

Lastly, we note that Algorithm \ref{Ho_QVI} requires setting a numerical tolerance which needs not be the same as the one of Algorithm \ref{A:Our_algorithm}. In fact, we have always taken it strictly smaller given the perceived unconditional convergence of the former and its precision.\footnote{Note that although precision is assessed in a different way in Algorithm \ref{Ho_QVI}, the pointwise residuals to the QVIs on the corresponding regions are afterwards checked as part of Algorithm \ref{A:Our_algorithm}}

\section{ Numerical results}
\label{S:Experiments}

In this section we assess the performance of Algorithm \ref{A:Our_algorithm}. We shall call the problems considered in the experiments simply `games'. The two value functions $V_1,V_2$ are approximated as described in Section \ref{section_our_algorithm} on an equispaced grid of $M+1$ nodes, and the computational domain is the plotted one. In all the subsequent games, $\varepsilon=10^{-8}$, $\alpha=0.8$ and $r_0=1$. The largest pointwise residual (equation (\ref{eq:R^max})) at convergence is denoted by $R^{\infty}$.

As mentioned before, the only NZSSIG for which an analytical solution is currently available is the benchmark game in Section \ref{Benchmark}. Therefore, we shall eventually focus on that problem. However, we introduce two other games first (for which we do not have an analytical solution). They illustrate how the initial guess for the benchmark game can be constructed and provide further numerical evidence supporting Algorithm \ref{A:Our_algorithm}.

\subsection{ Parabolic game}
This is a version of the benchmark game where the payoff is replaced by a concave parabola with roots $r_i^L$ and $r_i^R$:
\begin{equation}
\label{eq:parabolic_f}
\hat {f}_i(x)\defeq -(x-r_i^L)(x-r_i^R),\qquad r_i^L<r_i^R.
\end{equation}

Let us motivate this game. We seek to use as initial guess the value functions of the `unilateral games', that is, the control problems in which one of the players never intervenes (her continuation region is fixed and equal to $\mathbb R$). Removing the action of one of the players, however, may not always lead to a well-posed problem. Indeed, for payoffs without maximum one could end up for example with `infinite-valued value functions' or `infinite-valued optimal impulses'. This is indeed the case for the benchmark game. Thus, in order to skirt that difficulty it is convenient to define variations of the benchmark with payoffs that attain a maximum, like (\ref{eq:parabolic_f}). Note that, when well-posed, the unilateral games can be readily solved with Howard's algorithm (Algorithm \ref{Ho_QVI}).

\begin{figure}[H]
	\centering
	\includegraphics[scale=.4]{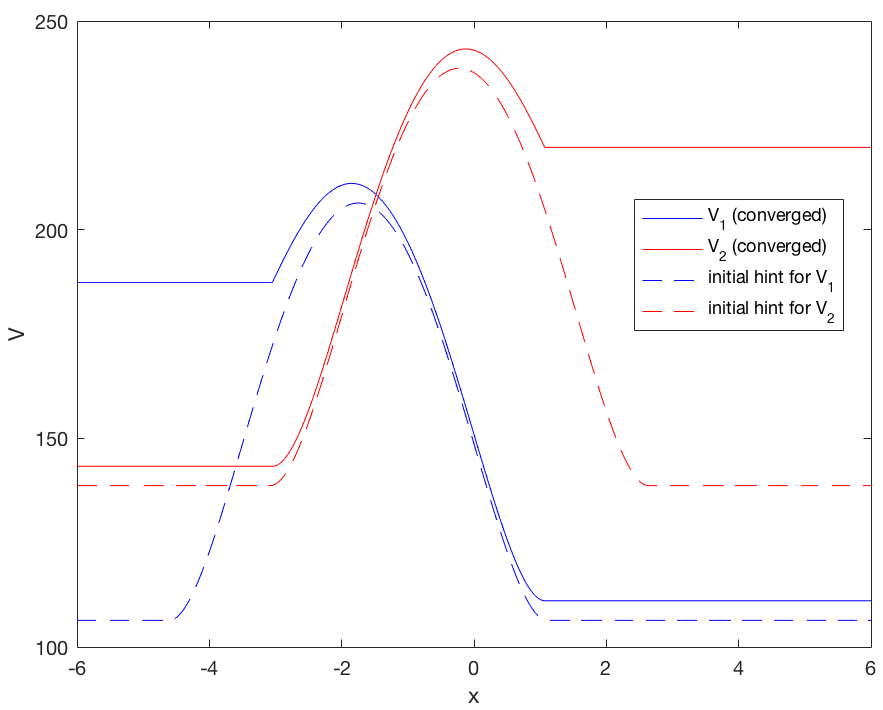}
	\caption{Value functions and Nash equilibrium with parabolic running payoffs: ${\hat f}_1=-(x+4.5)(x-1)$ and ${\hat f}_2=-(x+\pi)(x-2.7)$. Overlaid, initial guesses (solutions of the respective unilateral games)
		Parameters: $\rho=0.03$, $\sigma=0.25$, $c=100$, ${\tilde c}=30$, $\lambda={\tilde \lambda}=0$. 
		Here, $M=1000$; check Table \ref{T:Parabolic} for other values.}
	\label{F:Parabolic_V}
\end{figure}

Figure \ref{F:Parabolic_V} shows the numerical solution, $(V_1,V_2)$, to a parabolic game (pair of solid curves) along with the initial guess (pair of dashed curves). The latter are, in turn, the value functions of the corresponding unilateral games for players $1$ and $2$. It is intuitively clear that $(V_1,V_2)$ approximate, over the grid, functions which indeed satisfy the assumptions of the Verification Theorem \ref{verification}. The Nash equilibrium exhibited in this Theorem---$(\varphi_1^*,\varphi_2^*)$---can be retrieved from this graph. Indeed, note that in this case the cost for player $i$ is $\phi_i=-100$ and her approximate continuation region is $\{\mathbb M_i V_i-V_i<\varepsilon\}=\{\max V_i-V_i-100<\varepsilon\}$. 
Further, her optimal impulse at a given $x$ is the one that realizes (the discrete analogue of) the supremum in equation (\ref{Mi}), which in this case amounts to translating $x$ to the maximising point of $V_i$. Consequently, we get $\varphi_1^*=\big((-\infty,1.068), -1.848-x\big)$ and $\varphi_2^*=\big((-3.048,+\infty),-0.120-x\big)$. 

As it turns out, this parabolic game also converges from the `zero guess' ($V_1^0=V_2^0=0$). On Table \ref{T:Parabolic}, the convergence of $R^\infty$ over a wide range of finite difference grids (with $301,601,\ldots,3001$ nodes) is compiled. In all cases, both initial guesses lead to convergence (to within $\varepsilon$) of Algorithm \ref{A:Our_algorithm}. Nonetheless, the algorithm takes in general fewer iterations when starting from the values of the unilateral parabolic games. (We stress that those on Table \ref{T:Parabolic} are the outer iterations of Algorithm \ref{A:Our_algorithm}. Within each one there is an inner loop of Howard's iterations. Thus the computational cost scales linearly with the number of outer iterations.)

\begin{table}[H]
	\centering
	\begin{tabular}{llll}
		M & $R^\infty$ & its.$^a$ &  its.$^b$ \\ 
		\noalign{\smallskip}\hline\noalign{\smallskip}
		300  &  1.4$\times 10^{-12}$ & 54 & 53 \\
		600 &   4.0$\times 10^{-9}$ &  74 & 77 \\
		1200 & 3.3$\times 10^{-9}$ & 144 & 77 \\
		1800 & 9.7$\times 10^{-9}$& 95 & 77 \\
		2400 & 7.3$\times 10^{-9}$ & 123 & 77 \\
		3000 & 5.9$\times 10^{-9}$ & 215 & 103\\
		\noalign{\smallskip}\hline
	\end{tabular}
	\caption{Parabolic game: largest residual to QVIs at convergence ($R^\infty$) vs. grid points ($M+1$). Iterations to convergence within $\varepsilon=10^{-8}$ starting from: zero guess (its.$^a$) and value functions for one-player games (its.$^b$). (Same parameters as in Figure \ref{F:Parabolic_V}.)}
	\label{T:Parabolic} 
\end{table}

Since the exact solution to the parabolic game is unknown, we cannot say anything about the convergence of Algorithm \ref{A:Our_algorithm} in continuous-space. On the other hand, we see the system of QVIs is (approximately) enforced to within $\varepsilon\ll 1$ by a pair of numerical functions which---also approximately---comply with the regularity assumptions of Theorem \ref{verification}. As such, they are the (approximate) value functions of the parabolic game.

\begin{figure}[H]
	\centering
	\includegraphics[scale=.2]{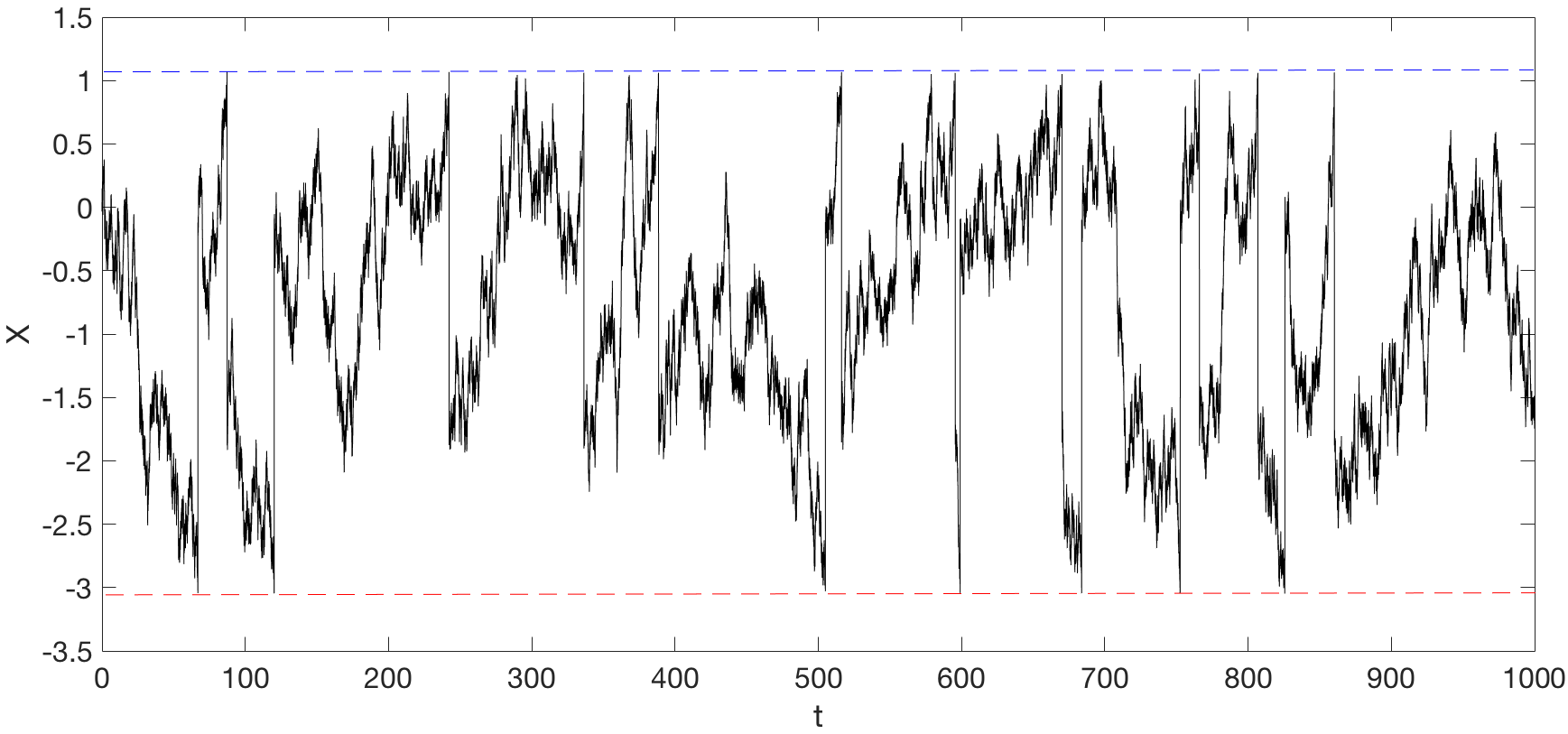}
	\caption{Exemplary trajectory from $x=0$ (parameters and solution from Figure \ref{F:Parabolic_V}). The blue and red dashed lines are the intervention thresholds for players 1 and 2, respectively.}
	\label{F:Trayectoria}
\end{figure}

We conclude this example by illustrating the interplay between the value functions which numerically solve the system of QVIs, the Nash equilibrium derived from them, and the evolution of the optimally controlled underlying. Once the optimal strategies $(\varphi_1^*,\varphi_2^*)$ are available, they can be executed on specific realizations of the game. Sticking to the parameters and numerical solution in Figure \ref{F:Parabolic_V}, Figure \ref{F:Trayectoria} depicts one exemplary trajectory of the underlying in the time interval $0\leq t\leq 1000$, starting from $x=0$ and subjected to the pair of optimal strategies. For numerical purposes, let us define
\begin{equation}
\label{J_aux}
\begin{split}
{\hat J}^i_T(x;\varphi_1,\varphi_2)& \defeq \mathbb E_x \left[ \int_0^T e^{-\rho_i s} f_i(X_s)ds + \sum_{k:\ \tau^i_k\leq T} e^{-\rho_i \tau^i_k} \phi_{i} \big( X_{(\tau^i_k)^-}, \delta^i_k \big) + \sum_{k:\ \tau^j_k\leq T} e^{-\rho_i \tau^j_k} \psi_{i} \big( X_{(\tau^j_k)^-}, \delta^j_k \big)\right].
\end{split}
\end{equation}

\begin{figure}[H]
	\centering
	\includegraphics[scale=.35]{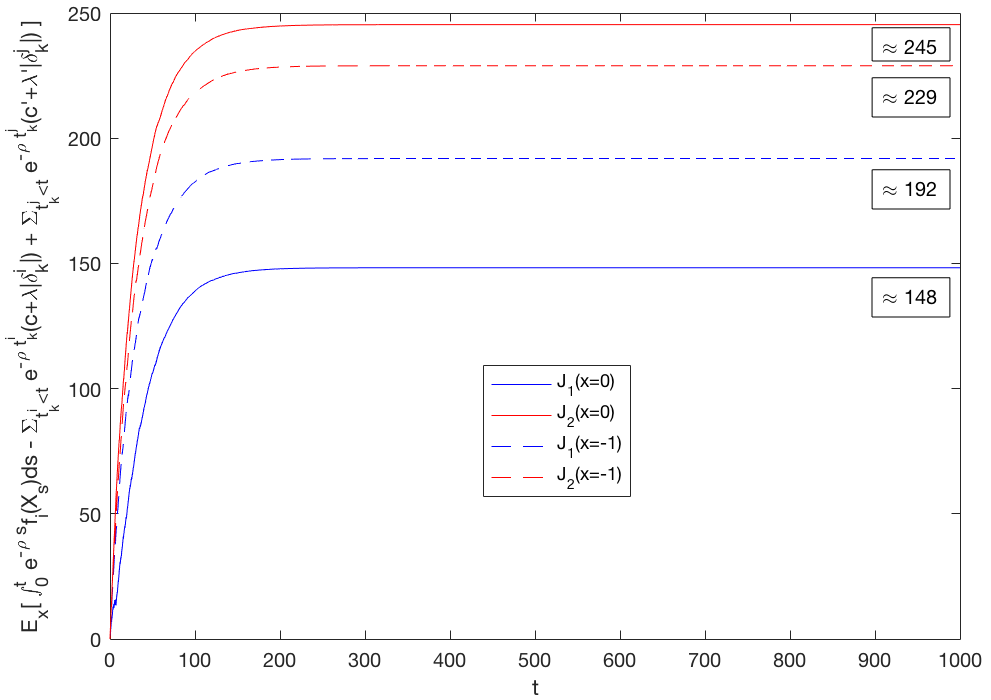}
	\caption{Approximations to the objective functions at $x=0$ (solid curves) and $x=-1$ (dashed curves), obtained by Monte Carlo simulation (see text for details). Parameters and optimal strategies from Figure \ref{F:Parabolic_V}. Compare with $V_1(0)$, $V_2(0)$, $V_1(-1)$ and $V_2(-1)$ there.}
	\label{F:Ex}
\end{figure}

\noindent Intuitively, ${\hat J}^i_T(x;\varphi_1,\varphi_2)\to J^i(x;\varphi_1,\varphi_2)$ (defined in (\ref{Vi_given_strats})) as $T\to\infty$. In fact, after $T\gtrsim 300$, the integrals in (\ref{J_aux}) for the parabolic game have essentially attained their asymptotic value. Thus, we simply take $J^i(x;\varphi_1,\varphi_2)\approx {\hat J}^i_{T=1000}(x;\varphi_1,\varphi_2)$. With this clarification, Figure \ref{F:Ex} shows in particular $V_1(0)=J^1(0;\varphi_1^*,\varphi_2^*)$, $V_2(0)=J^2(0;\varphi_1^*,\varphi_2^*)$,
	$V_1(-1)=J^1(-1;\varphi_1^*,\varphi_2^*)$, and
	$V_2(-1)=J^2(-1;\varphi_1^*,\varphi_2^*)$, obtained by Monte Carlo simulation.\footnote{The expected values in (\ref{J_aux}) are approximated by the mean over $N=200$ realizations integrated with the Euler-Maruyama method with time step $\Delta t=0.001$ } They compare fairly well with the values of $V_1(0)$, $V_2(0)$, $V_1(-1)$ and $V_2(-1)$ in Figure \ref{F:Parabolic_V} obtained by our algorithm. (Even better agreement could be obtained by increasing $M$ in that figure and reducing the discretisation bias and statistical error of the Monte Carlo simulation, but this is good enough to make our point.)

\begin{figure}[H]
	\centering
	\includegraphics[scale=.32]{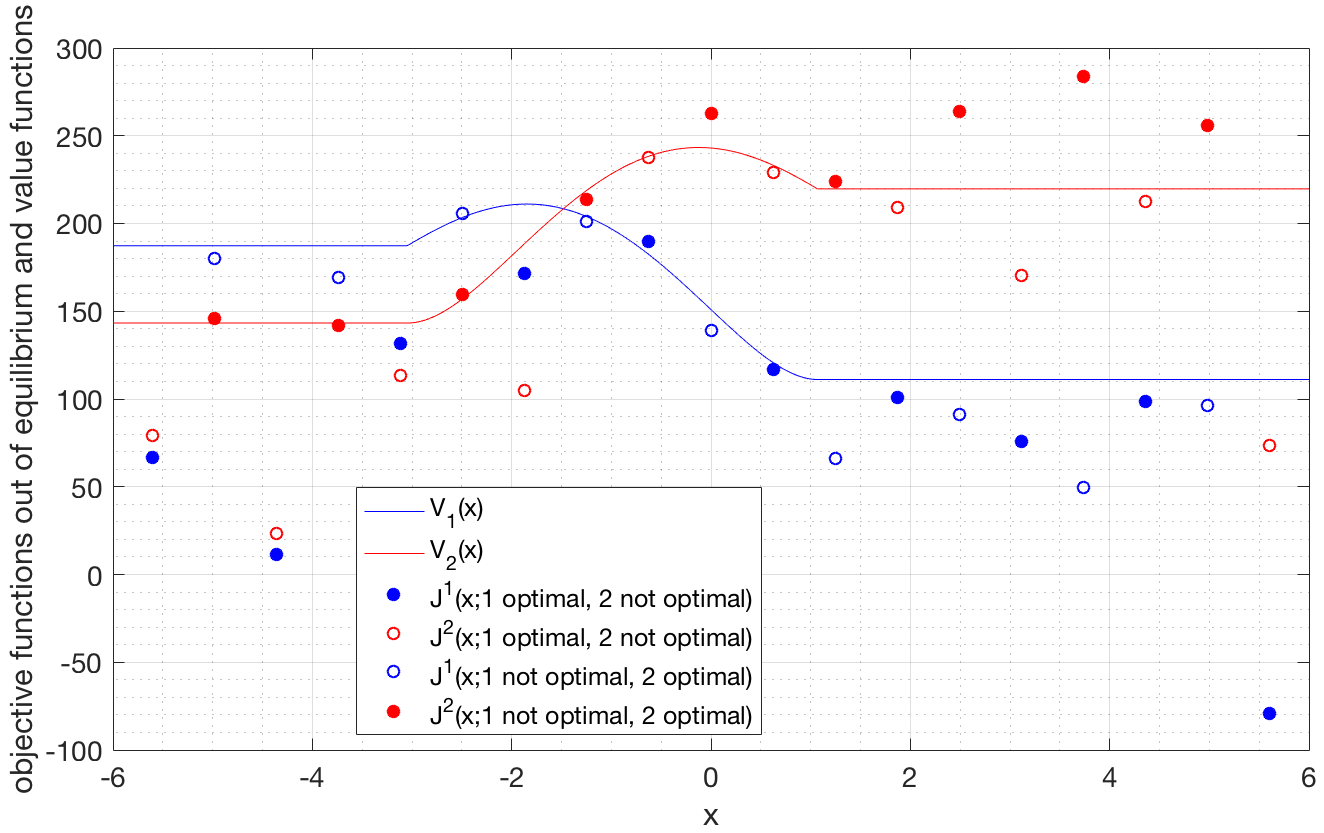}
	\caption{Empty blue [red] circles mean that the player 1 [2] has departed from her optimal strategy while her opponent has not. By virtue of Nash equilibrium, she cannot be (within numerical errors) better off than $V_1(x)$ [$V_2(x)$]. Full circles: objective function of the player who does not drop her optimal strategy. Parameters and optimal strategies from Figure \ref{F:Parabolic_V}.  (Note that results are subjected to numerical error.)}
	\label{F:Nash_eq}
\end{figure}

The Nash equilibrium itself can be visually explored in the following way. For a given starting point $x$, we keep the optimal strategy for one of the two players, and slightly alter the strategy of the other one. For concreteness, let us assume that player 1 "moves" (i.e. $\varphi_1=(1\pm 0.25{\cal U})\varphi^*_1$) while player 2 does not ($\varphi_2=\varphi^*_2$). (${\cal U}$ is the uniform distribution, and "$\pm$" means with equal chance.)
Then, we proceed to calculate $J^1(x;\varphi_1,\varphi_2^*)$ by Monte Carlo simulation as before. By definition of the Nash equilibrium, $J^1(x;\varphi_1,\varphi_2^*)$ can not be larger than $V_1(x)$. Within numerical tolerance, this is indeed observed in Figure \ref{F:Nash_eq}, where the blue empty circles (representing $J^1(x;\varphi_1,\varphi_2^*)$) do not lie over the blue solid curve (which represents $V_1(x)$). Full blue circles represent $J^2(x;\varphi_1,\varphi_2^*)$: note that the player who sticks to her optimal strategy may indeed improve over $V_2(x)$, should her opponent depart from a Nash equilibrium. (When player 2 is the one who changes, the red circles and red curve in Figure \ref{F:Nash_eq} apply instead.)

We stress, however, that Monte Carlo simulations such as those cannot {\em prove} that a pair of strategies form a Nash equilibrium. (At most, they could disprove it.) The only way---in the current state of the theory---of fully characterizing a Nash equilibrium calls for solving the system of QVIs---which Algorithm \ref{A:Our_algorithm} has now made possible.

\setcounter{subsection}{1}
\subsection{ Capped benchmark game}

In this game, we replace the running payoffs of the benchmark by a version capped at $K>0$:
\begin{equation}
\label{eq:capped_f}
{\bar f}_i(x)\defeq\min\{(-1)^{i-1}(x-s_i), K\}
\end{equation}
We shall always take $K=5$. Once again, the corresponding unilateral games are well-posed and their solutions can be used as initial guess for the capped benchmark game. As in the previous example, the capped game also seems to converge from the zero guess. Some convergence results are compiled in Table \ref{T:Capped}. Note that convergence falters with $M=600$ and $M=3300$ (independently of $k_{\max}$); we will come back to this later.

\begin{table}[H]
	\centering
	\begin{tabular}{lll}
		M & $R^\infty$ & its.\\
		
		\noalign{\smallskip}\hline\noalign{\smallskip}
		600 &---& $\infty$ \\
		900 & $6.5\times 10^{-10}$ & 100 \\
		1200 & $8.8\times 10^{-9}$ & 124 \\
		1500 & $7.4\times 10^{-9}$ & 78 \\
		2700 & $4.1\times 10^{-9}$ & 96 \\
		3000 & $6.7\times 10^{-9}$ & 125 \\
		3300 &---& $\infty$\\
		\noalign{\smallskip}\hline
	\end{tabular}
	\caption{Convergence of $R^\infty$ in the capped benchmark game using the zero guess. The hyphen stands for lack of convergence within $\varepsilon$. Same parameters as in Figure \ref{F:Capped_error}.}
	\label{T:Capped} 
\end{table}

\begin{figure}[H]
	\centering
	\includegraphics[scale=.33]{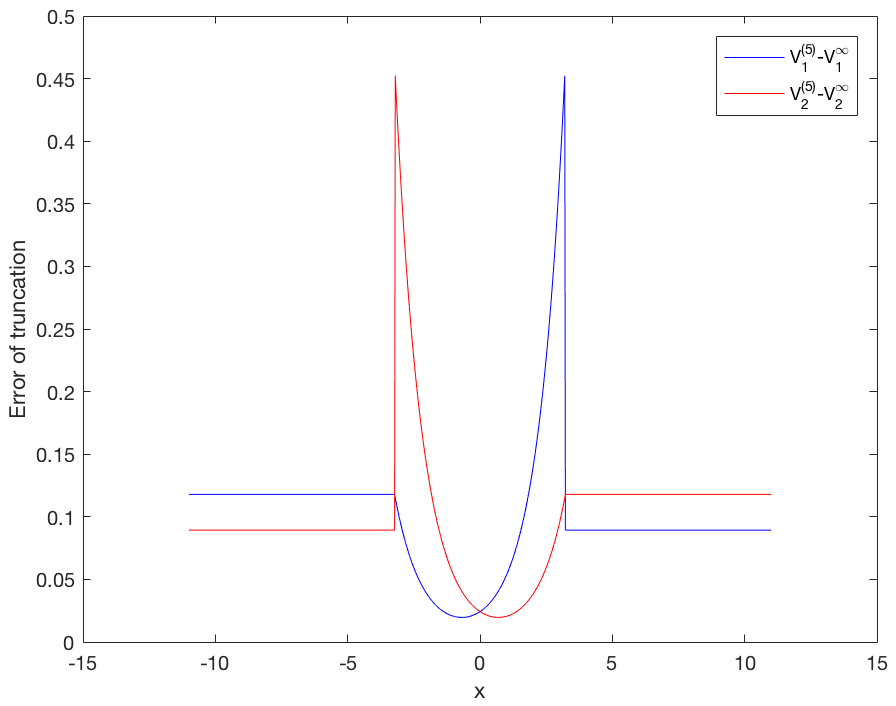}
	\caption{Difference (in absolute value, less than $1\%$) between the (numerical, with $M=1200$) value functions of the capped benchmark game and those (exact) of the benchmark game---justifying the former being used as initial guess for the latter. Parameters are: $\sigma=0.25, \rho=0.03, c=100, {\tilde c}=30, \lambda=0.5, {\tilde\lambda}=0.3,  s_1=-\pi/3, s_2=\pi/3$ and $K=5$ (for capped game).}
	\label{F:Capped_error}
\end{figure}

The value functions of the capped benchmark game are a good approximation to those of the benchmark game itself (see Figure \ref{F:Capped_error}). This seems to make sense: due to the action of the opponent and for $\sigma\ll K$, the discarded portion of the payoff is not very relevant in practice.

\subsection{ Benchmark game with educated initial guess}
Finally, we tackle the benchmark game, for which an exact solution is available (see Section \ref{Benchmark}). Contrary to the previous examples, Algorithm \ref{A:Our_algorithm} does not seem to enjoy unconditional convergence here. In fact, when the zero guess was used, it failed to converge more often than not (not reported). In order to construct an adequate initial guess, we first solve for the value functions of the capped benchmark game. 
Using them as the initial guess, convergence was achieved in every experiment.

\begin{figure}[H]
	\hspace*{-.85cm}
	\includegraphics[scale=.33]{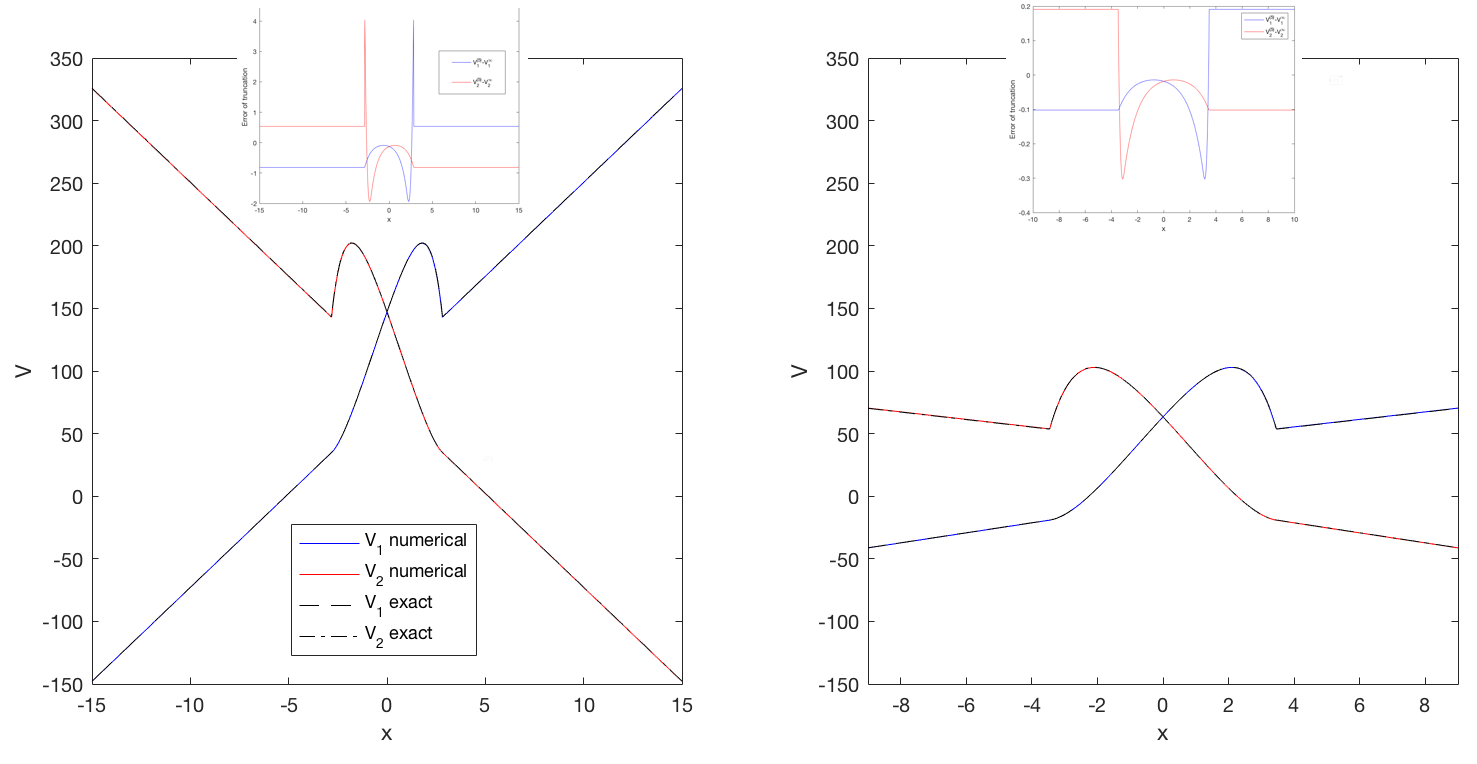}
	\caption{Value functions and Nash equilibrium of two instances of benchmark game. Initial guess: solutions of capped benchmark games ($K=5$). Overlaid, error of the initial guess.
		Parameters: $\rho=0.02$, $\sigma=0.15$, $s_1=-3$, $s_2=3$, $c=100$, ${\tilde c}=0$, $\lambda={\tilde \lambda}=15$,  
		$M=1000$  (left); $\rho=0.03$, $\sigma=0.25$, $s_1=-2$, $s_2=2$, $c=100$, ${\tilde c}=30$, $\lambda=4$, ${\tilde \lambda}=3$, 
		$M=1000$  (right).} 
	\label{F:Benchmark_V}
\end{figure}

The result of two experiments are plotted on Figure \ref{F:Benchmark_V}. The numerical approximations can hardly be distinguished from the exact solutions with the naked eye. As it was done with the parabolic game, once again we can retrieve an approximate Nash equilibrium from the numerical solution. For the left figure, this is $\varphi_1^*=\big((-2.82,+\infty),1.53-x\big)$ and $\varphi_2^*=\big((-\infty,2.82),-1.53-x\big)$. When compared with the exact equilibrium, the errors on the corresponding abscissae are smaller than $0.15$, i.e., smaller than half the grid step.

\begin{table}[H]
	\centering
	\begin{tabular}{llll}
		M & $R^\infty$ & $|$error$|$ & its.\\
		\noalign{\smallskip}\hline\noalign{\smallskip}
		500 & 3.8$\times 10^{-10}$ & 0.687 & 126 \\
		1000 & 1.2$\times 10^{-9}$ & 0.805 & 154 \\
		1500 & 7.2$\times 10^{-9}$ & 0.512 & 157 \\
		2000 & 3.5$\times 10^{-9}$ & 0.365 & 172 \\
		2500 &---&---& $\infty$\\ 
		\noalign{\smallskip}\hline
	\end{tabular}
	\caption{Convergence of Algorithm \ref{A:Our_algorithm} for benchmark game. Same parameters as in Figure \ref{F:Benchmark_V} (left).}
	\label{T:Benchmark_ejemploB} 
\end{table}

\begin{table}[H]
	\centering
	\begin{tabular}{llll}
		M & $R^\infty$ & $|$error$|$ & its.\\
		\noalign{\smallskip}\hline\noalign{\smallskip}
		600 &---&---& $\infty$ \\
		800 & 9.5$\times 10^{-10}$ & 0.023 & 183 \\
		1000 & 3.7$\times 10^{-9}$ & 0.330 & 177 \\
		1400 & 8.7$\times 10^{-9}$ & 0.196 & 159 \\
		1800 & 6.8$\times 10^{-9}$ & 0.121 & 226 \\
		2200 & 5.6$\times 10^{-9}$ & 0.073 & 224 \\
		2600 &---&---& $\infty$\\ 
		\noalign{\smallskip}\hline
	\end{tabular}
	\caption{Convergence of Algorithm \ref{A:Our_algorithm} for benchmark game (same parameters as in Figure \ref{F:Benchmark_V} (right).}
	\label{T:Benchmark_ejemploD} 
\end{table}

\begin{figure}[H]
\centering
	\includegraphics[scale=.6]{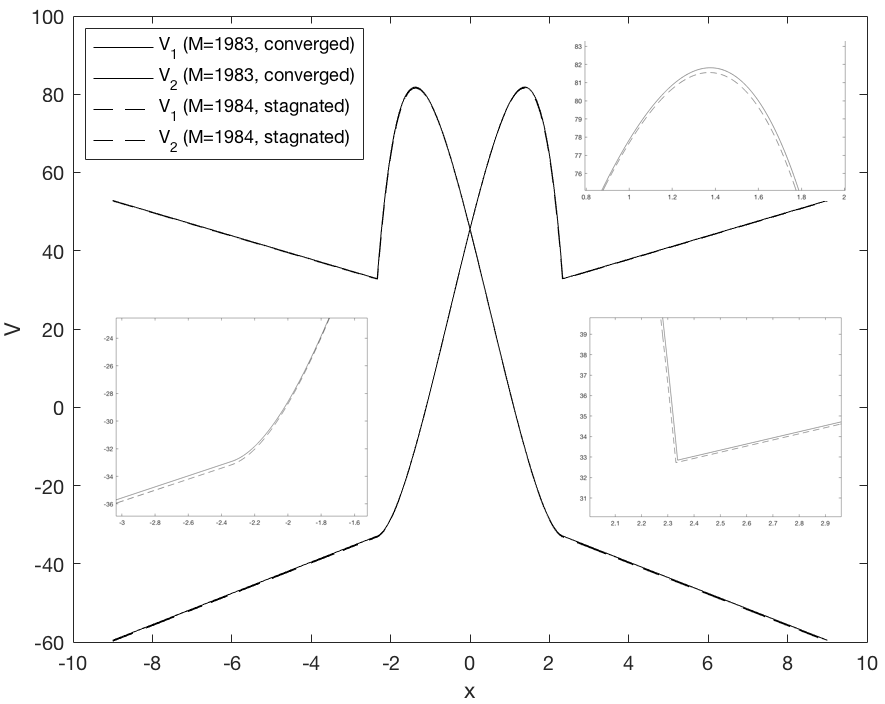}
	\caption{Stagnation of accuracy: the pair with $M=1983$ is fully converged (within numerical tolerance); the pair with $M=1984$ (overlaid) is not. The insets zoom in on both pairs of functions. Parameters: $\rho=0.02$, $\sigma=0.15$, $s_1=-1$, $s_2=1$, $c=100$, ${\tilde c}=30$, $\lambda=4$, ${\tilde\lambda}=3$.} 
	\label{F:Stagnation}
\end{figure}

On Tables \ref{T:Benchmark_ejemploB} and \ref{T:Benchmark_ejemploD}, the convergence of the numerical approximation provided by Algorithm~ \ref{A:Our_algorithm} to the true solution is demonstrated. However, $R^\infty$ fails to drop below $\varepsilon$ for a fine enough discretisation. This reflects a pattern: $R^\infty$ stagnates as $M$ grows. Upon closer inspection, it turns out that the stagnating largest pointwise residual for each player takes place at the junction between the intervention and continuation regions of the opponent, where the exact value function of the former player has a singularity (a \textit{kink}). As discussed in Section \ref{section_our_algorithm}, Algorithm \ref{A:Our_algorithm} is always going to place those kinks at finite difference nodes. The overshooting residuals are thus due to the inability of a numerical solution to reproduce a sharp, nonsmooth feature---where the finite difference derivatives grow unbounded as the distance between nodes goes to zero. Figure \ref{F:Stagnation} illustrates this situation. It can be seen that errors continue to be acceptable (far less than $1\%$ relative error in the worst case). We also highlight the fact that the kinks do not seem to bring about oscillations of the numerical solutions around them---the notorious Gibbs' phenomenon.

\section{ Conclusions}
We have designed and tested a novel policy iteration algorithm---the first one as far as we know---to numerically solve nonzero-sum stochastic impulse games (NZSSIGs). The approach consists in solving a system of quasi-variational inequalities which characterizes the value functions and Nash equilibrium, exploiting a recent
theoretical breakthrough in \cite{ABCCV_V2}.

Our algorithm computes iteratively the approximate
solution by partitioning, for each of the players, the discretised spatial domain into an
approximate continuation region and an approximate
intervention region. They are defined through a relaxation
parameter that evolves along the iterations. In the  continuation
region, we solve one quasi-variational inequality by means of (a generalization of) Howard's algorithm, whereas in the complement, a gain is computed. A strategy for producing an educated initial guess to start the iterations---which relies on solving two associated, standard impulse control problems---has been presented along with the new Algorithm as well. 
 
We have not carried out a convergence analysis. Instead, we have gathered plenty of numerical evidence showing that the Algorithm can be applied
confidently. In particular, we have  performed convergence tests on the only available---at the time of writing---example of an analytically solvable NZSSIG. In other test NZSSIGs for which we do not have an analytic solution,
we have explored consistency and numerical results on the value functions and Nash equilibrium, also with satisfactory results.

Value functions of NZSSIGs can develop sharp kinks at the confluence between the continuation and the intervention regions of the opponent. Capturing such features into numerical approximations is a pervasive challenge of numerical analysis. Numerical tests show that the presence of such kinks may eventually put a cap to the accuracy attained by our algorithm. On the other hand, its stability is not affected by them. Moreover, in every case the largest pointwise error was perfectly acceptable for the purposes of most applications.

In sum, the new Algorithm offers a means of gaining quantitative insight into applications modelled by NZSSIGs. Natural continuations of the present work include the convergence analysis of the Algorithm,
and enriching the discretisation method so as to better capture singularities in the solutions.

\section*{ Acknowledgments}
JPZ was supported by CNPq, FAPERJ, and the Brazilian-French network in Mathematics.

\noindent FB gratefully acknowledges support from the  Finance for Energy Market Research Centre (FiME).

\bibliographystyle{amsalpha}
\bibliography{references}
\addcontentsline{toc}{section}{References}

\end{document}